\documentclass[a4paper,11pt]{article}

\usepackage[latin1]{inputenc}
\usepackage{latexsym}

\usepackage{amssymb}
\usepackage{amsfonts}
\usepackage{amsmath}
\usepackage{amsthm}
\usepackage{textcomp}

\usepackage{comment}

\title{Tropicalization of group representations}
\author{Daniele Alessandrini \\ \small \ \\ \small \textit{Universit\`a di Pisa, Italy} \\ \small \textit{E-mail address:} daniele.alessandrini@gmail.com}
\date{}

\usepackage[dvips]{graphicx}
\newcommand{\figuresize}{4.9cm}

\newtheorem{theorem}{Theorem}[section]    
\newtheorem{lemma}[theorem]{Lemma}          
\newtheorem{corollary}[theorem]{Corollary}
\newtheorem{proposition}[theorem]{Proposition}

\theoremstyle{definition}
\newtheorem{definition}[theorem]{Definition}    

\newcommand{\nuovo}[1]{{{\bfseries \upshape #1}}}

\newcommand{\erre}{{\mathbb{R}}}
\newcommand{\ci}{{\mathbb{C}}}

\newcommand{\pro}{{\mathbb{P}}}
\newcommand{\cappa}{{\mathbb{K}}}
\newcommand{\effe}{{\mathbb{F}}}

\newcommand{\iper}{{\mathbb{H}}}
\newcommand{\ident}{{\mbox{Id}}}
\newcommand{\tro}{{\mathbb{T}}}

\newcommand{\ocors}{{\mathcal{O}}}

\newcommand{\gfamil}{{\mathcal{G}}}

\newcommand{\ecors}{{\mathcal{E}}}

\newcommand{\fcors}{{\mathcal{F}}}

\DeclareMathOperator{\trace}{tr}

\DeclareMathOperator{\charat}{Char}

\DeclareMathOperator{\Span}{Span}
\DeclareMathOperator{\Imm}{Im}

\newenvironment{gmatrice}{ \begin{pmatrix} }{ \end{pmatrix} }
\newenvironment{pmatrice}{ \left( \begin{smallmatrix} }{ \end{smallmatrix}  \right) }

\newcommand{\pvettore}[1]{{ \begin{pmatrice} #1 \end{pmatrice} }}

\newcommand{\freccia}{{\longrightarrow}}

\begin{document}

\sloppy

\maketitle

\begin{abstract}    
In this paper we give an interpretation to the boundary points of the compactification of the parameter space of convex projective structures on an $n$-manifold $M$. These spaces are closed semi-algebraic subsets of the variety of characters of representations of $\pi_1(M)$ in $SL_{n+1}(\erre)$. The boundary was constructed as the ``tropicalization'' of this semi-algebraic set. Here we show that the geometric interpretation for the points of the boundary can be constructed searching for a tropical analogue to an action of $\pi_1(M)$ on a projective space. To do this we need to construct a tropical projective space with many invertible projective maps. We achieve this using a generalization of the Bruhat-Tits buildings for $SL_{n+1}$ to non-archimedean fields with real surjective valuation. In the case $n=1$ these objects are the real trees used by Morgan and Shalen to describe the boundary points for the Teichm\"uller spaces. In the general case they are contractible metric spaces with a structure of tropical projective spaces.
\end{abstract}


\section{Introduction}

Let $M$ be a closed oriented $n$-manifold such that $\pi_1(M)$ is virtually centerless, it is Gromov-hyperbolic and it is torsion-free. We denote by $\mathcal{T}^{c}_{\erre\pro^n}(M)$ the parameter space of marked convex projective structures on $M$. If $S$ is an orientable hyperbolic surface of finite type, we denote by $\mathcal{T}^{cf}_{\iper^2}(S)$ the Teichm\"uller space of $S$.

In \cite{A2} we showed that the space $\mathcal{T}^{c}_{\erre\pro^n}(M)$ can be identified with a closed semi-algebraic subset of the character variety $\charat(\pi_1(M),SL_{n+1}(\erre))$. Then we applied the Maslov dequantization to this semi-algebraic set (see also \cite{A1}) and, using an inverse limit of logarithmic limit sets of this space, we constructed the tropical counterpart of $\mathcal{T}^{c}_{\erre\pro^n}(M)$. The spherical quotient of this tropical counterpart, denoted by $\partial \mathcal{T}^{c}_{\erre\pro^n}(M)$, can be glued to $\mathcal{T}^{c}_{\erre\pro^n}(M)$ ``at infinity'', defining a compactification $\mathcal{T}^{c}_{\erre\pro^n}(M) \cup \partial \mathcal{T}^{c}_{\erre\pro^n}(M)$ of the parameter space. The same construction applied to the Teichm\"uller space $\mathcal{T}^{cf}_{\iper^2}(S)$ gives back the Thurston boundary $\partial \mathcal{T}^{cf}_{\iper^2}(S)$.

The aim of this paper is to give a geometric interpretation of the points of these tropical counterparts. We are guided by the idea that the points of the tropicalization of a parameter space should be related with the tropical counterparts of the parametrized objects. Every point of $\mathcal{T}^{c}_{\erre\pro^n}(M)$ corresponds to a conjugacy class of representations of $\pi_1(M)$ in $SL_{n+1}(\erre)$. Geometrically such a representation corresponds to an action of $\pi_1(M)$ on a vector space $\erre^{n+1}$, or, equivalently, on a projective space $\erre\pro^n$. In this paper we introduce the tropical counterparts of these actions, i.e. actions of a group on tropical modules and tropical projective spaces. This is the notion we propose of tropicalization of a group representation.

There is a naif notion of tropical projective space, the projective quotient of a free module $\tro^n$, but these spaces have few invertible projective maps, hence they have few group actions. We give a more general notion of tropical modules and, correspondingly, of tropical projective spaces. We show that these objects have an intrinsic metric, the tropical version of the Hilbert metric, that is invariant for tropical projective maps, and that the topology induced by this metric is contractible. Then we construct a special class of tropical projective spaces, denoted by $P^n$, by using a generalization of the Bruhat-Tits buildings for $SL_{n+1}$ to non-archimedean fields with a surjective real valuation.

In the usual case of a field $\effe$ with discrete valuation, Bruhat and Tits constructed a polyhedral complex of dimension $n$ with an action of $SL_{n+1}(\effe)$. In the case $n=1$, Morgan and Shalen generalized this construction to a field with a general valuation, and they studied these objects using the theory of real trees. We extend this to general $n$, and we think that the proper structure to study these objects is the structure of tropical projective spaces. The paper \cite{JSY}, developed independently from this work, contains a similar approach to the Bruhat-Tits buildings. Tropical geometry is used there to study the convexity properties of the Bruhat-Tits buildings for $SL_n(\effe)$, for a field $\effe$ with discrete valuation.

With every point of the boundary we can associate a class of representations of $\pi_1(M)$ in $SL_{n+1}(\effe)$, where $\effe$ is real closed non-archimedean field with a surjective real valuation. Every representation of $\pi_1(M)$ in $SL_{n+1}(\effe)$ induces an action by tropical projective maps on our tropical projective spaces $P^n$. We compute the length spectra of these actions on $P^n$, and we show that the length spectrum of an action identifies a boundary point in $\partial \mathcal{T}^{c}_{\erre\pro^n}(M)$. Then we use the fact that tropical projective spaces are contractible to show that for every action of the fundamental group of the manifold on a tropical projective space there exists an equivariant map from the universal covering of the manifold to the tropical projective space.

This theorem can hopefully lead to interesting consequences about the interpretation of the boundary points. For example in the case $n=1$, where $P^1$ is a real tree, the equivariant map induces a duality between actions of the fundamental group on $P^1$ and measured laminations on the surface (see \cite{MS84}, \cite{MS88} and \cite{MS88'}). It would be very interesting to extend this result to the general case. For example an action of the fundamental group of the surface on a tropical projective space $P^n$ induces a degenerate metric on the surface, and this metric can be used to associate a length with each curve. Anyway it is not clear up to now how to classify these induced structures. This is closely related to a problem raised by J. Roberts (see \cite[problem 12.19]{Oh01}): how to extend the theory of measured laminations to higher rank groups, such as, for example, $SL_n(\erre)$.

A brief description of the following sections. In section \ref{sec:first definitions} we give elementary definitions of semifields, semimodules and projective spaces over a semifield, and we give some examples of semimodules.

In section \ref{sec:Linear maps between free semimodules} we discuss invertibility of linear maps in $\tro^n$ and the tropicalizations of linear maps on a vector space $\effe^n$ over a non archimedean field $\effe$. With every such map $f$ we associate a linear map $f^\tau$ on $\tro^n$, and we discuss the relations between $f^{\tau}$ and ${(f^{-1})}^{\tau}$: globally they are not inverse one of the other, but this happens on a specific ``inversion-domain''.

In section \ref{sec:Tropical projective structure on Bruhat-Tits buildings} we define the structure of tropical projective space we put on the generalization of the Bruhat-Tits buildings, and we give a description of this space. Tropical modules $\tro^n$ can be seen as the tropicalization of a vector space $\effe^n$ over a non-archimedean field $\effe$, but this tropicalization depends on the choice of a basis of $\effe^n$. Our description with tropical charts, one for each basis of $\effe^n$, can be interpreted by thinking the Bruhat-Tits buildings as a tropicalization of $\effe^n$ with reference to all possible bases.

In section \ref{sec:metric spaces} we define in a canonical way a metric on tropical projective spaces making tropical segments geodesics and tropical projective maps $1$-Lipschitz. This metric is the transposition to tropical geometry of the Hilbert metric on convex subsets of $\erre\pro^n$. The topology induced by this metric is shown to be contractible.

Finally, in section \ref{sec:Tropicalization of group actions} we consider a representation of a group $\Gamma$ in $SL_{n+1}(\effe)$, and we study the induced action by tropical projective maps on our Bruhat-Tits building. First we compute the length spectrum of the action with reference to the canonical metric, and, if $\Gamma = \pi_1(M)$, we show show how we can recover the information characterizing a boundary point. Then by using the fact that tropical projective spaces are contractible, we show that every action of $\pi_1(M)$ on a tropical projective space has an equivariant map from the universal cover of $M$ to the tropical projective space.

\section{First definitions}     \label{sec:first definitions}

\subsection{Tropical semifields}

We need some linear algebra over the tropical semifield. By a \nuovo{semifield} we mean a quintuple $(S,\oplus,\odot,0,1)$, where $S$ is a set, $\oplus$ and $\odot$ are associative and commutative operations $S\times S \freccia S$ satisfying the distributivity law, $0, 1\in S$ are, respectively, the neutral elements for $\oplus$ and $\odot$. Moreover we require that every element of $S^* = S\setminus \{0\}$ has a multiplicative inverse. We will denote the inverse of $a$ by $a^{\odot-1}$. Given an element $b \neq 0$ we can write $a\oslash b = a \odot b^{\odot-1}$. Note that $0$ is never invertible and $\forall s \in S: 0\odot s = 0$.

A semifield is called \nuovo{idempotent} if $\forall s \in S: s\oplus s = s$. In this case a partial order relation is defined by 
$$ a \leq b \Leftrightarrow a\oplus b = b$$
We will restrict our attention to the idempotent semifields such that this partial order is total. In this case $(S\setminus \{0\},\odot,\leq)$ is an abelian ordered group.
Vice versa, given an abelian ordered group $(\Lambda,+,<)$, we add to it an extra element $-\infty$ with the property $\forall \lambda \in \Lambda: -\infty < \lambda$, and we define a semifield:
$$\tro = \tro_\Lambda = (\Lambda\cup\{-\infty\}, \oplus, \odot, -\infty, 0)$$
with the tropical operations $\oplus, \odot$ defined as 
$$a\oplus b = \max(a,b)$$ 
$$a\odot b = \left\{ 
\begin{array}{ll}
a+b & \mbox{ if } a,b\in \Lambda\\
-\infty & \mbox{ if } a=-\infty \mbox{ or } b=-\infty
\end{array}
\right. $$

We will use the notation $1_\tro = 0$, as the zero of the ordered group is the one of the semifield, and $0_\tro = -\infty$. If $a\in \tro$ and $a \neq 0_\tro$, then $a \odot (-a) = 1_\tro$.
Hence $-a = a^{\odot-1}$, the \nuovo{tropical inverse} of $a$. The order on $\Lambda\cup\{-\infty\}$ induces a topology on $\tro$ that makes the operations continuous.

Semifields of the form $\tro = \tro_\Lambda$ will be called \nuovo{tropical semifields}. The semifield that in literature is called the tropical semifield is, in our notation, $\tro_\erre$.

We are interested in tropical semifields because they are the images of valuations. Let $\effe$ be a field, $\Lambda$ an ordered group, and $v:\effe \freccia \Lambda \cup \{+\infty\}$ a surjective valuation.
Instead of using the valuation, we prefer the \nuovo{tropicalization map}: 
$$\tau:\effe \ni z \freccia -v(z) \in \tro = \tro_\Lambda = \Lambda \cup \{-\infty\}$$ 
The tropicalization map satisfies the properties of a norm:
\begin{enumerate}
\item $\tau(z) = 0_\tro \Leftrightarrow z = 0$
\item $\tau(zw) = \tau(z)\odot \tau(w)$
\item $\tau(z+w) \leq \tau(z) \oplus \tau(w)$
\item $\tau$ is surjective.
\end{enumerate}

For every element $\lambda \in \tro$ we choose an element $t_\lambda \in \effe$ such that $\tau(t_\lambda) = \lambda$.

We will denote the valuation ring by $\ocors=\{z\in\effe \ |\ \tau(z) \leq 1_\tro\}$, its unique maximal ideal by $m=\{z\in\effe \ |\ \tau(z) < 1_\tro\}$, its residue field by $D=\ocors / m$ and the projection by $\pi:\ocors\freccia D$.

\begin{proposition} The map $\tau$ `often' sends $+$ to $\oplus$, i.e.:
\begin{enumerate}
\item $\tau(w_1) \neq \tau(w_2) \Rightarrow \tau(w_1+w_2)=\tau(w_1)\oplus\tau(w_2)$. 

\item If $\tau(w_1)=\tau(w_2)=\lambda$, then $t_{-\lambda} w_1, t_{-\lambda} w_2 \in \ocors \setminus m$. In this case\\ $\pi(t_{-\lambda}w_1) + \pi(t_{-\lambda}w_2) \neq 0 \in D \Rightarrow \tau(w_1+w_2)=\tau(w_1)\oplus\tau(w_2)$. 
\end{enumerate}
\end{proposition}

\begin{proof}
It follows from elementary properties of valuations.
\end{proof}

\subsection{Tropical semimodules and projective spaces}

\begin{definition}
Given a semifield $S$, an $S$-\nuovo{semimodule} is a triple $(M,\oplus,\odot,0)$, where $M$ is a set, $\oplus$ and $\odot$ are operations:
$$\oplus:M\times M \freccia M \hspace{2cm} \odot:S\times M \freccia M$$ 
$\oplus$ is associative and commutative and $\odot$ satisfies the usual associative and distributive properties of the product by a scalar. We will also require that
$$\forall v \in M: 1\odot v = v  \hspace{2cm}  \forall v \in M: 0 \odot v = 0$$
\end{definition}

Note that the following properties also holds:
$$\forall a \in S: a\odot 0 = 0$$ 
$$\forall a \in S^*: \forall v\in M: a\odot v = 0 \Rightarrow v = 0$$ 
The first follows as $a\odot 0 \oplus b = a\odot 0 \oplus a\odot(a^{-1}\odot b) = a \odot (0 \oplus a^{-1}\odot b) = a \odot a^{-1}\odot b = b$. And then the second follows as $0=a\odot v \Rightarrow 0 = a^{-1}\odot 0 = 1 \odot v = v $.

Most definitions of linear algebra can be given as usual. Let $S$ be a semifield and $M$ a $S$-semimodule. A \nuovo{submodule} of $M$ is a subset closed for the operations. If $v_1, \dots, v_n \in M$, a \nuovo{linear combination} of them is an element of the form $c_1 \odot v_1 \oplus \dots \oplus c_n \odot v_n$. If $A\subset M$ is a set, it is possible to define its \nuovo{spanned submodule} $\Span_S(A)$ as the smallest submodule containing $A$ or, equivalently, as the set of all linear combinations of elements in $A$.  A \nuovo{linear map} between two $S$-semimodules is a map preserving the operations. The image of a linear map is a submodule, but (in general) there is not a good notion of kernel.

If $S$ is an idempotent semifield, then $M$ is an idempotent semigroup for $\oplus$. In this case a partial order relation is defined by 
$$ v \leq w \Leftrightarrow v\oplus w = w$$
Linear maps are monotone with reference to this order.

Let $S$ be a semifield and $M$ be an $S$-module. The \nuovo{projective equivalence} relation on $M$ is defined as:
$$x \sim y \Leftrightarrow \exists \lambda \in S^* : x = \lambda\odot y$$
This is an equivalence relation. The \nuovo{projective space} associated with $M$ may be defined as the quotient by this relation:
$$\pro(M) = (M \setminus \{0\}) / \sim $$
The quotient map will be denoted by $\pi:M\setminus \{0\} \freccia \pro(M)$.

The image by $\pi$ of a submodule is a \nuovo{projective subspace}.

If $f:M\freccia N$ is a linear map, we note that $v \sim w \Rightarrow f(v) \sim f(w)$. The linear map induces a map between the associated projective spaces provided that the following condition holds:
$$\{v\in M \ |\ f(v)=0 \} \subset \{0\}$$
We will denote the induced map as $\overline{f}: \pro(M) \freccia \pro(N)$. Maps of this kind will be called \nuovo{projective maps}. The condition does not imply in general that the map is injective. Actually a projective map $\overline{f}: \pro(M) \freccia \pro(M)$ may be not injective nor surjective in general.

\subsection{Examples}

From now on we will consider only semimodules over a tropical semifield $\tro = \tro_\Lambda$. The simplest example of $\tro$-semimodule is the \nuovo{free} $\tro$-semimodule of \nuovo{rank} $n$, i.e. the set $\tro^n$ where the semigroup operation is the component wise sum, and the product by a scalar is applied to every component. If $x \in \tro^n$ we will write by $x^1, \dots, x^n$ its components:
$$ x = \pvettore{x^1\\ \vdots\\ x^n}$$

These modules inherit a topology from the order topology of the tropical semifields: the product topology on the free modules and the subspace topology on their submodules. The partial order on these semimodules can be expressed in coordinates as
$$\forall x,y \in \tro^n: x \preceq y \Leftrightarrow \forall i: x^i \leq y^i$$

Other examples are the submodules 
$$F\tro^n = \Span_\tro({(\tro^*)}^n ) = {(\tro^*)}^n \cup \{0_\tro\} \subset \tro^n$$

The projective space associated with $\tro^n$ is $\pro(\tro^n) = \tro\pro^{n-1}$, and the projective space associated with $F\tro^n$ is $\pro(F\tro^n) = F\tro\pro^{n-1}$. 
We will denote its points with homogeneous coordinates: $$\pi(x) = [x^1:x^2:\dots:x^n]$$

These projective spaces inherit the quotient topology, and projective maps are continuous for this topology.

$\tro\pro^1 = \pro(\tro^2)$ can be identified with $\Lambda \cup \{-\infty, +\infty\}$ via the map:
$$\tro\pro^1 \ni [x^1:x^2] \freccia x^1 - x^2 \in \Lambda \cup \{-\infty, +\infty\}$$
With this identification $\tro\pro^1$ inherits an order: given $a = [a^1:a^2], b = [b^1:b^2] \in \tro\pro^1$, we define $a \preceq b \Leftrightarrow a^1-a^2 \leq b^1 - b^2$. All tropical projective maps $\tro\pro^1 \freccia \tro\pro^1$ are never increasing or never decreasing with reference to this order.
We give a name to three special points: $0_\tro = [0_\tro:1_\tro] = -\infty, 1_\tro = [1_\tro:1_\tro] = 0, \infty_\tro = [1_\tro:0_\tro] = +\infty$.

When $\Lambda=\erre$, $\tro_\erre\pro^{n-1}$ may be described as an $(n-1)$-simplex, whose set of vertices is $\{\pi(e_1),\dots,\pi(e_n)\}$ ($e_i$ being the elements of the canonical basis of $\tro^n$). Given a set of vertices $A$, the face with vertices in $A$ is the projective subspace $\pi(\Span_\tro(A))$. $F\tro\pro^{n-1}$ is naturally identified with the interior of the simplex $\tro\pro^{n-1}$.

\section{Linear maps between free semimodules}    \label{sec:Linear maps between free semimodules}

\subsection{Tropical matrices}

As before let $\tro = \tro_\Lambda$ be a tropical semifield. Let $e_i$ be the element of $\tro^n$ having $1$ as the $i$-th coefficient and $0$ as the others. These elements form the \nuovo{canonical basis} of $\tro^n$.

Let $f:\tro^n \freccia \tro^m$ be a linear map. Then we can define the matrix $A=[f]=(a^i_j)$ as $a^i_j = {(f(e_j))}^i$. The usual properties of matrices and linear maps hold in this case:
\begin{enumerate}
\item $f(e_j)$ is the $j$-th column of $[f]$: $f(e_j)= \bigoplus_i a^i_j \odot e_i$.
\item If $v\in \tro^n$, ${(f(v))}^i = \bigoplus_j a^i_j \odot v^j$, or $f(v)= \bigoplus_{i,j} a^i_j \odot v^j\odot e_i$.
\item $f$ is surjective $\Leftrightarrow$ the columns of $[f]$ span $\tro^m$.
\item There is a binary correspondence between linear maps and matrices with entries in $S$.
\item The matrix of the composition of two maps is the product matrix, i.e. $[f \circ g]=[f]\odot[g]$, where ${(A\odot B)}^i_j = \bigoplus_k A^i_k \odot B^k_j$. 
\end{enumerate}

The identity matrix, corresponding to the identity map $\ident_\tro:\tro^n\freccia\tro^n$, will be also denoted by $\ident_\tro = ({(\delta_\tro)}^i_j)$, where 
$${(\delta_\tro)}^i_j = 
\left\{ \mbox{
\begin{tabular}{ll}
$1_\tro$ &  if $i=j$\\ 
$0_\tro$ & if $i\neq j$
\end{tabular}
}\right.$$

A linear map $f:\tro^n \freccia \tro^m$ induces a linear map $f:F\tro^n \freccia F\tro^m$ by restriction, provided that no element in $F\tro^n$ is mapped outside $F\tro^m$, i.e. if every row of the matrix $[f]$ contains a non-zero element.

Projective maps $\overline{f}:\tro\pro^{n-1} \freccia \tro\pro^{m-1}$ are induced by matrices mapping no non-zero vector to zero. These are precisely the matrices such that every column contains a non-zero element.

Tropical linear maps are very seldom surjective. This depends on the following property:
$$\Span_\tro(v_1,\dots,v_m) = \tro^m \Leftrightarrow \forall i=1,\dots,m: \exists a \in \tro^*: a \odot e_i \in \{v_1, \dots, v_m\}$$
    
Hence a tropical linear map is surjective if and only if it has, among its columns, all the elements of the canonical basis of the codomain.

Let $f:\tro^n \freccia \tro^m$ be a linear map, with matrix $[f]=(a^i_j)$. Suppose that every column of $[f]$ contains a non-zero element. We will denote by $f^{pi}:\tro^m \freccia \tro^n$ the map defined by:
$${(f^{pi}(y))}^j = \min_i (y^i-a^i_j)$$
(in the previous formula, by $-0_\tro$ we mean an element greater than every other element in $\tro$. This value is never the minimum, thanks to the condition on the columns). In \cite{CGQ04} this map is called residuated map.

\begin{theorem}       \label{teo:pseudoinverse}
Let $y\in \tro^m$. Then $y \in \Imm f$ if and only if exists a sequence $\epsilon:\{1,\dots,m\}\freccia \{1,\dots,n\}$ such that
$$\forall k=1,\dots,m : y^k-a^k_{\epsilon_k} ={(f^{pi}(y))}^{\epsilon_k}$$
Moreover we have
$$f^{-1}(y) = \bigcup_{\epsilon \mbox{\scriptsize\  as before}} \left\{ x \in \tro^n \ |\ 
\mbox{
\begin{tabular}{l}
$x \preceq f^{pi}(y) $ \\
$\forall k=1,\dots,m : x^{\epsilon_k}= {(f^{pi}(y))}^{\epsilon_k}$  
\end{tabular}
}
\right\}
$$
This implies that $f^{-1}(y)$ is a single point if and only if every function $\epsilon$ as before is surjective.

The function $f^{pi}$ plays the role of a pseudo-inverse function, as it sends every point of the image in one of its pre-images, in a continuous way. It has the following properties:
\begin{enumerate}
\item $\forall x\in \tro^n: \forall y \in \tro^m: \left(x \preceq f^{pi}(y) \Leftrightarrow f(x) \preceq  y\right)$.
\item $\forall x \in \tro^n: x \preceq f^{pi}(f(x))$
\item $\forall y \in \tro^m: f(f^{pi}(y)) \preceq y$
\item $\forall y \in \Imm f: f(f^{pi}(y)) = y$
\item $\forall x \in \Imm f^{pi}: f^{pi}(f(x)) = x$
\item $f_{|\Imm f^{pi}}:\Imm f^{pi} \freccia \Imm f$ is bijective, with inverse $f^{pi}$.
\end{enumerate}
\end{theorem}

\begin{proof}
The point $y$ is in the image if and only if exists $x\in \tro^n$ such that $f(x)=y$. Then
$$ f(x)=y \ \Leftrightarrow\  \forall i : \bigoplus_j (a^i_j \odot x^j) = y^i \ \Leftrightarrow\ 
\left\{
\mbox{
\begin{tabular}{l}
$\forall i,j: a^i_j + x^j \leq y^i $\\
$\forall i: \exists j: a^i_j + x^j = y^i$  
\end{tabular}
}
\right.
\ \Leftrightarrow $$
$$\left\{
\mbox{
\begin{tabular}{l}
$\forall i,j: x^j \leq y^i-a^i_j $\\
$\forall i: \exists j: x^j = y^i - a^i_j$  
\end{tabular}
}
\right.
\ \Leftrightarrow\ 
\left\{
\mbox{
\begin{tabular}{l}
$\forall j: x^j \leq \displaystyle\min_i (y^i-a^i_j) $\\
$\forall i: \exists j: x^j = y^i - a^i_j$  
\end{tabular}
}
\right. $$
Then
$$ y \in \Imm f \ \Leftrightarrow \exists \epsilon : \forall k : y^k-a^k_{\epsilon_k} = \min_i (y^i-a^i_{\epsilon_k})$$
In this case $x^{\epsilon_k} = y^k-a^k_{\epsilon_k}$. All the claims of the theorem follows from the calculations above.
\end{proof}

\subsection{Simple tropicalization of linear maps}

Let $\effe$ be a valued field, with tropicalization map $\tau:\effe \freccia \tro$. An $\effe$-vector space $\effe^n$ may be tropicalized through the componentwise tropicalization map, again denoted by $\tau:\effe^n \freccia \tro^n$.

Let $f:\effe^n \freccia \effe^m$ be a linear map, expressed by a matrix $[f]=(a^i_j)$. Its tropicalization is the map $f^{\tau}:\tro^n \freccia \tro^m$ defined by the matrix $[f^{\tau}]=(\alpha^i_j)=(\tau(a^i_j))$.

\begin{proposition} The following properties hold:
\begin{enumerate}
\item $\forall z \in \effe^n: \tau(f(z)) \preceq f^\tau(\tau(z))$.
\item $\forall x \in \tro^n: \exists z \in \effe^n: \tau(z)=x$ and $\tau(f(z)) = f^\tau(x)$.
\end{enumerate}
\end{proposition}

Let $A \in GL_n(\effe)$ be an invertible matrix. Its tropicalization $\alpha=A^\tau:\tro^n \freccia \tro^n$ (i.e. $\alpha=(\alpha^i_j)=(\tau(a^i_j))$) is, in general, not invertible. Anyway it has the property that every column and every row contains a non-zero element, hence it has a pseudo-inverse function, and it induces a linear map $F\tro^n\freccia F\tro^n$, and projective maps $\tro\pro^{n-1}\freccia\tro\pro^{n-1}$ and $F\tro\pro^{n-1}\freccia F\tro\pro^{n-1}$.

Now let $B=A^{-1}$, the inverse of $A$. We will write $\beta=B^\tau$. We would like to see $\beta$ as an inverse of $\alpha$, but this is impossible, as $\alpha$ is not always invertible. 

\begin{proposition} The following statements hold
\begin{enumerate}
\item $\forall i,j: {(\alpha \odot \beta)}^i_j \geq {(\delta_\tro)}^i_j$ and ${(\beta\odot\alpha)}^i_j \geq {(\delta_\tro)}^i_j$.
\item $\forall x \in \tro^n: x \preceq \alpha(\beta(x))$ and $y \preceq \beta(\alpha(y))$.
\item $\forall x \in \tro^n: \alpha^{pi}(x) \preceq \beta(x)$.
\item $\forall x \in \tro^n : \alpha(\beta(x)) = x \Leftrightarrow \beta(x) = \alpha^{pi}(x) $ 
\end{enumerate}
\end{proposition}

\begin{proof}
\begin{enumerate}
\item It follows from: $AB=\ident$, $BA=\ident$, $\tau(z_1+z_2) \leq \tau(z_1)\oplus\tau(z_2)$.
\item It follows from the previous statement.
\item This is equivalent to $\forall i: \max_j (\beta^i_j + x^j) \geq \min_j (x_j - \alpha^j_i)$, i.e. $\forall i: \exists k,h: \beta^i_k + x^k \geq x^h - \alpha^h_i$. This always holds as, from the first statement, we know that $\max_k (\beta^i_k+\alpha^k_i) = {(\beta \odot \alpha)}^i_i \geq 1_\tro$, hence $\forall i: \exists k: \beta^i_k + x^k \geq x^k - \alpha^k_i$.  
\item From theorem \ref{teo:pseudoinverse}, part 1, we know that $\alpha(\beta(x)) \preceq x \Leftrightarrow \beta(x) \preceq \alpha^{pi}(x)$. The reversed inequalities always holds.
\end{enumerate}
\end{proof}

If $\alpha$ and $\beta$ are tropicalizations of two maps $A,B \in GL_n(\effe)$ such that $A^{-1} = B$, we will call \nuovo{inversion domain} the set $D_{\alpha\beta} = \{x \in \tro^n \ |\ \alpha(\beta(x))=x\}$.

\begin{proposition} \label{prop:inversion domains}
The inversion domains have this name because of the following property: $D_{\beta\alpha}=\beta(D_{\alpha\beta})$, $D_{\alpha\beta} = \alpha(D_\beta\alpha)$ and $\beta_{|D_{\alpha\beta}}:D_{\alpha\beta}\freccia D_{\beta\alpha}$ is bijective with inverse $\alpha_{|D_{\beta\alpha}}:D_{\beta\alpha} \freccia D_{\alpha\beta}$.

The set $D_{\alpha\beta}$ is a tropical submodule, and we can write explicit equations for it:
$$D_{\alpha\beta} = \{ x \in \tro^n \ |\ \forall h,k: x^h-x^k \geq {(\alpha\odot\beta)}^h_k \}$$
As a consequence if $A \in GL_n(\ocors)$, then $D_{\alpha\beta}\neq \emptyset$. Note that the matrices $\alpha$ and $\beta$ are not one the inverse of the other, but, in the hypothesis $D_{\alpha\beta}\neq \emptyset$, then $\forall i:{(\alpha\odot\beta)}^i_i = 1_\tro$.

The map $\beta_{|D_{\alpha\beta}}$ is the composition of a permutation of coordinates and a tropical dilatation: there exists a diagonal matrix $d$ and a permutation of coordinates $\sigma$ such that ${(\sigma \circ d \circ \beta)}_{|D_{\alpha\beta}} = \ident_{|D_{\alpha\beta}}:D_{\alpha\beta}\freccia D_{\alpha\beta}$.
\end{proposition}

\section{Tropical projective structure on Bruhat-Tits buildings}   \label{sec:Tropical projective structure on Bruhat-Tits buildings}

\subsection{Definition}

Given a non-archimedean field $\effe$ with a surjective real valuation, we are going to construct a family of tropical projective spaces we will call $P^{n-1}(\effe)$, or simply $P^{n-1}$ when the field is well understood. This family arises as a generalization of the Bruhat-Tits buildings for $SL_n$ to non-archimedean fields with surjective real valuation. In the usual case of a field with integral valuation, Bruhat and Tits constructed a polyhedral complex of dimension $n-1$ with an action of $SL_n(\effe)$. In the case $n=2$, Morgan and Shalen generalized this construction to a field with a general valuation, and they studied these objects using the theory of real trees. We want to extend this to general $n$, and we think that the proper structure to study these objects is the structure of tropical projective spaces.

Let $V = \effe^n$, an $\effe$-vector space of dimension $n$ and an infinitely generated $\ocors$-module. We consider the natural action $GL_n(\effe) \times V \freccia V$.

\begin{definition} 
An $\ocors$-\nuovo{lattice} of $V$ is an $\ocors$-finitely generated $\ocors$-submodule of $V$. 
\end{definition}

\begin{proposition}
Let $L$ be an $\ocors$-finitely generated $\ocors$-submodule of $V$. Then every minimal set of generators is $\effe$-linearly independent, hence $L$ is free. 
\end{proposition}

\begin{proof}
Let $\{e_1, \dots, e_m\}$ be a minimal set of generators of $L$. Suppose they are not $\effe$-independent. Then there exist $a_1, \dots, a_m \in \effe$ s.t. $\sum a_i e_i = 0$. We may suppose  $\tau(a_1) \leq \dots \leq \tau(a_m)$. There exist elements $b_1, \dots, b_m\in \ocors$ s.t. $a_i = b_i a_m$. Hence $a_m ( \sum b_i e_i ) = 0 \Rightarrow e_m = \sum_1^n b_i e_i$ with $b_1, \dots, b_{m-1} \in \ocors$. They can't be minimal.  

An element of $L$ is an $\ocors$-linear combination of $\{e_1, \dots, e_m\}$ because they are generators, and the linear combination is unique because they are $\effe$-independent. Hence $L$ is free.
\end{proof}

If $L$ is a finitely generated $\ocors$-submodule of $V$, its rank is a number from $0$ to $n$. 

\begin{definition}
A \nuovo{maximal} $\ocors$-\nuovo{lattice} is an $\ocors$-lattice of rank $n$. 
\end{definition}

We denote by $U^n(\effe)$ (or simply $U^n$) the set of all $\ocors$-lattices of $V=\effe^n$, and by $FU^n(\effe)$ (or simply $FU^n$) the subset of all maximal $\ocors$-lattices and the $\ocors$-lattice $\{0\}$.

$U^n$ and $FU^n$ can be turned in $\tro$-semimodules by means of the following operations:

$\oplus : U^n \times U^n \freccia U^n \hspace*{2cm} L \oplus M = \Span_\ocors(L \cup M)$

$\odot : \tro \times U^n \freccia U^n \hspace*{2.25cm}  x \odot L = z L\mbox{, where }z \in \effe, \tau(z) = x$

The associated tropical projective spaces will be denoted by $\pro(U^n(\effe)) = P^{n-1}(\effe)$ and $\pro(FU^n(\effe)) = FP^{n-1}(\effe)$. We will simply write $P^{n-1}$ and $FP^{n-1}$ when the field $\effe$ is understood.

As we said there is a natural action $GL_n(\effe)\times V \freccia V$. Every element $A\in GL_n(\effe)$ sends $\ocors$-lattices in $\ocors$-lattices, hence we have an induced action $GL_n(\effe)\times U^n \freccia U^n$. This action preserves the rank of a lattice, and in particular it sends $FU^n$ in itself. Among the $\ocors$-lattices with the same rank this action is transitive, for example there exist an $A \in GL_n(\effe)$ sending every maximal $\ocors$-lattice of $V$ in the standard lattice $\ocors^n \subset V$.  

Hence the group $SL_n(\effe)$ acts naturally on $U^n$ and $FU^n$ by tropical linear maps and on $P^{n-1}$ and $FP^{n-1}$ by tropical projective maps.

\subsection{Description}

Let $\ecors=(e_1, \dots, e_n)$ be a basis of $V$. We denote by $\varphi_\ecors:\tro^n \freccia U^n$ the map:
$$\varphi_\ecors(y)=\varphi_\ecors(y^1,\dots,y^n) = I_{y^1}e_1 + \dots + I_{y^n}e_n = \Span_\ocors(t_{y^1}e_1, \dots t_{y^n}e_n)$$

\begin{proposition}       \label{corol:olattices}
Let $<e_1, \dots, e_m>$ be a $\ocors$-basis of an $\ocors$-lattice $L$, and let $p_i \in \effe$. Then:
\begin{enumerate}
	\item $p_i e_i \in L$ $\Leftrightarrow$ $p_i \in \ocors$ $\Leftrightarrow$ $\tau(p_i) \leq 0$.  
	\item $<p_1 e_1, \dots p_n e_n>$ is an $\ocors$-basis of $L$ $\Leftrightarrow$ $p_i \in \ocors \setminus m$ $\Leftrightarrow$ $\tau(p_i) = 0$.
\end{enumerate}  
\end{proposition}

\begin{proof}
It follows from the properties of valuations.
\end{proof}

This proposition implies that $\varphi_\ecors$ is injective and $\varphi_\ecors(F\tro^n) \subset FU^n$. For every basis $\ecors$ we have a different map $\varphi_\ecors$. The union of the images of all these maps is the whole $U^n$, and the union of all the sets $\varphi(F\tro^n)$ is equal to $FU^n$. We will call the maps $\varphi_\ecors$ \nuovo{tropical charts} for $U^n$. Theorem \ref{teo:charts} will justify this name. 

Note that the charts respect the partial order relations on $\tro^n$ and on $U^n$:
$$x \preceq y \Leftrightarrow \varphi(x) \subset \varphi(y)$$

\begin{lemma}  \label{lemma:common basis}
Let $L,M \subset V$ be two $\ocors$-lattices, and suppose that $L$ is maximal. Then there is a basis $v_1,\dots,v_n$ of $L$ and scalars $a_1,\dots,a_n \in \effe$ such that $a_1 v_1, \dots a_n v_n$ is a basis of $M$.
\end{lemma}

\begin{proof}
Fix a basis $e_1,\dots,e_n$ of $V$ such that $L = \Span_\ocors(e_1,\dots,e_n)$. Let $f_1,\dots, f_n$ be a basis of $M$. For every vector $f_i$ there is a scalar $b_i\in \effe$ such that $b_i f_i \in L$. Then, if $b_i$ is the one with maximal valuation, $b_i M \subset L$. The thesis follows by applying \cite[corol. II.3.2]{MS84} to the $\ocors$-modules $L$ and $b_i M$.
\end{proof} 

\begin{corollary}        \label{corol:common chart}
Given two points $x,y \in U^n$, there is a tropical chart containing both of them in its image. 
\end{corollary}

Given two bases $\ecors=(e_1,\dots,e_n)$ and $\fcors=(f_1,\dots,f_n)$, we have two charts $\varphi_\ecors, \varphi_\fcors$. We want to study the intersection of the images.

We put $I=\varphi_\ecors(\tro^n) \cap \varphi_\fcors(\tro^n)$, $I_\ecors = \varphi_\ecors^{-1}(I)$, $I_\fcors = \varphi_\fcors^{-1}(I)$. We want to describe the sets $I_\fcors,I_\ecors$ and the \nuovo{transition function}: $\varphi_{\fcors\ecors}=\varphi_\fcors^{-1} \circ \varphi_\ecors: I_\ecors \freccia I_\fcors$.

The transition matrices between $\ecors$ and $\fcors$ are denoted by $A=(a^i_j),B=(b^i_j)\in GL_n(\effe)$: 
$$\forall j: e_j=\sum_i a^i_j f_i \ \ \ \ \ \ \ \ \forall j: f_j=\sum_i b^i_j e_i \ \ \ \ \ \ \ \ A=B^{-1}$$

We will write $\alpha=A^\tau$ and $\beta=B^\tau$, i.e. $\alpha=(\alpha^i_j)=(\tau(a^i_j))$, $\beta=(\beta^i_j)=(\tau(b^i_j))$. 

\begin{theorem}[{[Description of the tropical charts]}]    \label{teo:charts}
We have that $I_\fcors = D_{\alpha\beta}$ and $I_\ecors = D_{\beta\alpha}$, the inversion domains described in proposition \ref{prop:inversion domains}. Moreover $\varphi_{\fcors\ecors}=\alpha_{|I_\ecors}$ and $\varphi_{\ecors\fcors}=\beta_{|I_\fcors}$, the tropicalizations of the transition matrices.
\end{theorem}

\begin{proof}
First, we need to prove the following two assertions:
\begin{enumerate}
\item $\varphi_\ecors(y)\subset \varphi_\fcors(x) \Leftrightarrow \alpha(y) \preceq x$ and $\varphi_\fcors(x)\subset \varphi_\ecors(y) \Leftrightarrow \beta(x) \preceq y$.
\item $\varphi_\fcors(x)=\varphi_\ecors(y) \Leftrightarrow  x=\alpha(y) \mbox{ and } y=\beta(x) $.
\end{enumerate}

Let $w = \sum_i w^i f_i \in \effe^n$. Then:

$w\in\Span_\ocors(f_1,\dots,f_n) \Leftrightarrow \forall i: w^i \in \ocors \Leftrightarrow \forall i: \tau(w^i)\leq 1_\tro$

$t_y w\in \varphi_\fcors(x) \Leftrightarrow \forall i: \frac{t_y w^i}{t_{x^i}}\in\ocors \Leftrightarrow  \forall i: \tau(w^i)\leq x^i-y$

$\varphi_\ecors(y) \subset \varphi_\fcors(x) \Leftrightarrow \forall j:t_{y^j}e_j \in \varphi_\fcors(x)\Leftrightarrow \forall j,i: \tau(a^i_j)\leq x^i-y^j$

$\varphi_\fcors(x) \subset \varphi_\ecors(y) \Leftrightarrow  \forall j,i: \tau(b^i_j)\leq y^i-x^j$

Then we have: 
$$\varphi_\ecors(y) = \varphi_\fcors(x) \Leftrightarrow \forall j,i: \tau(b^i_j)\leq y^i-x^j \leq -\tau(a^j_i)$$ $$\Leftrightarrow \forall j,i: \tau(b^i_j)+x^j \leq y^i \leq -\tau(a^j_i)+x^j$$
$$ \Leftrightarrow \forall i: \max_j(\tau(b^i_j)+x^j) \leq y^i \leq \min_j(-\tau(a^j_i)+x^j)$$ 
$$\Leftrightarrow \forall i: \bigoplus_j (\beta^i_j \odot x^j) \leq y^i \leq \min_j(x^j-\alpha^j_i) \Leftrightarrow \forall i: {(\beta(x))}^i \leq y^i \leq {(\alpha^{pi}(x))}^i$$ 

The map $\varphi_\ecors$ is injective, hence, given a fixed $x$, if an $y$ satisfying the last condition exists, it has to be unique. Then the interval in which its coordinates are free to vary must degenerate to a single point. Then we have:

$$\varphi_\ecors(y) = \varphi_\fcors(x) \Leftrightarrow \beta(x)=y=\alpha^{pi}(x)$$ 

We can prove the symmetric equalities reversing the roles of $\ecors$ and $\fcors$.

Now we look at $\varphi_\fcors^{-1}(I) = \{x\in\tro^n \ |\ \exists y \in \tro^n: \varphi_\ecors(y) = \varphi_\fcors(x)\}$. We have 
$$x \in \varphi_\fcors^{-1}(I) \Leftrightarrow \exists y: \forall i: \max_j(\beta^i_j+x^j) \leq y^i \leq \min_j(-\alpha^j_i+x^j)$$
$$\Leftrightarrow \forall i: \max_j(\beta^i_j+x^j) \leq \min_j(-\alpha^j_i+x^j) \Leftrightarrow \forall i,k,h: \beta^i_k+x^k \leq -\alpha^h_i+x^h$$
$$\Leftrightarrow \forall i,k,h: x^h-x^k\geq \beta^i_k+\alpha^h_i \Leftrightarrow \forall k,h: x^h-x^k \geq \oplus_i(\alpha^h_i \odot \beta^i_k)$$
\end{proof}

\section{Tropical projective spaces as metric spaces}   \label{sec:metric spaces}

\subsection{Finitely generated semimodules}

Free semimodules have the usual universal property: let $M$ be a $\tro$-semimodule, and $v_1, \dots, v_n \in M$. Then there is a linear map:

\begin{center}
\vspace{0.2cm}
\begin{tabular}{ccc}
$\tro^n$ & $\freccia$ & $\Span_\tro(v_1,\dots,v_n)$\\
 \\
$c$ & $\freccia$ & $c^1 \odot v_1 \oplus \dots \oplus c^n \odot v_n$ 
\end{tabular}
\vspace{0.2cm}
\end{center}

This map sends $e_i$ in $v_i$ and its image is $\Span_\tro(v_1,\dots,v_n)$.

Hence every finitely generated $\tro$-semimodule is the image of a free $\tro$-semimodule.

In the following we will need some properties of finitely generated semimodules over $\tro_\erre$, so for this section we will suppose $\tro = \tro_\erre$.

First we want to discuss a pathological example we prefer to neglect. Consider the following equivalence relation on $\tro^2$:

$$(x^1,x^2) \sim (y^1,y^2)  \Leftrightarrow \left\{ 
\begin{array}{l}
x^1 < x^2, y^1 < y^2 \mbox{ and } x^2 = y^2\\
or\\
x^1 \geq x^2, y^1 \geq y^2 \mbox{ and } x^1 = y^1
\end{array}
\right.$$

\begin{figure}[htbp]
\begin{center}
\includegraphics[height=\figuresize]{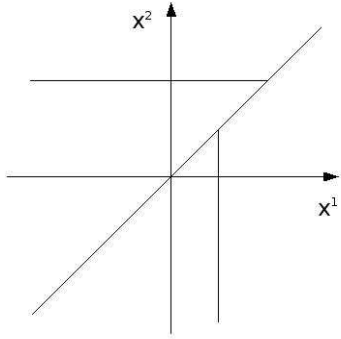}
\caption{Two examples of equivalence classes for the relation defining the quotient module $B$: $\{ x^2 = 2, x^1 < 2\}$ and $\{ x^1 = 1, x^2 \leq x^1\}$.}  \label{fig:equiv_class}
\end{center}
\end{figure}

The quotient for this relation will be denoted by $B$.
If $a\sim a'$ and $b\sim b'$, then $a\oplus b = a' \oplus b'$ and $\lambda \odot a = \lambda \odot a'$. Hence the operations $\oplus$, $\odot$ induces operations on $B$, turning it in a finitely generated $\tro$-semimodule. We will denote the equivalence classes in the following way: if $(x^1,x^2)$ satisfies $x^1 < x^2$ we will denote its class as $[(\cdot,x^2)]$, if $x^1 \geq x^2$ we will denote its class as $[(x^1,\cdot)]$. The $\odot$ operation act as $$\lambda\odot [(\cdot,x^2)] = [(\cdot,\lambda\odot x^2)]$$
and analogously for the other classes. The $\oplus$ operation acts as 
$$[(\cdot,x^2)] \oplus [(\cdot,y^2)] = [(\cdot,x^2 \oplus y^2)]$$
$$[(x^1,\cdot)] \oplus [(y^1,\cdot)] = [(x^1 \oplus y^1),\cdot]$$
$$[(x^1,\cdot)] \oplus [(\cdot,x^2)] = [(x^1,x^2)]$$

If we put on the quotient a topology making the projection continuous, then the point $[(x^1,\cdot)]$ is not closed, as its closure must contain the point $[(\cdot,x^1)]$.

We define a $\tro$-semimodule to be \nuovo{separated} if it does not contain any submodule isomorphic to $B$. We will see in the following section that every separated $\tro$-semimodule has a natural metrizable topology making all linear maps continuous. Examples of separated $\tro$-semimodules are all free semimodules (as there exists no submodule in $\tro^n$ whose associated  projective space has exactly two points) and the semimodules $U^n$ (as every two points in $U^n$ are in the image of the same tropical chart, hence in a submodule isomorphic to $\tro^n$).

\begin{lemma}
Let $M$ be a $\tro$-semimodule and let $f:\tro^2 \freccia M$ be a linear map such that $f\pvettore{x^1\\x^2} = f\pvettore{y^1\\x^2} = m$ when $y^1 \leq x^1$. Then $\forall y \leq x^1: f\pvettore{y\\x^2} = m$.
\end{lemma}

\begin{proof}
Case 1): If $y^1 \leq y \leq x^1$, then $\pvettore{y^1\\x^2} \preceq \pvettore{y\\x^2} \preceq \pvettore{x^1\\x^2}$. Linear maps are monotone with reference to $\preceq$, hence $m \preceq f\pvettore{y\\x^2} \preceq m$.

Case 2): If $y = y^1 - (x^1-y^1)$, then consider the points $a = (y^1-x^1) \odot \pvettore{x^1\\x^2}$ and $b = (y^1-x^1) \odot\pvettore{y^1\\x^2}$. We have $f(a) = f(b) = (y^1-x^1) \odot m$. Then 
$$f\pvettore{y\\x^2} = f\left(\pvettore{y\\x^2} \oplus b\right) = f\pvettore{y\\x^2} \oplus f(b) = $$ 
$$ = f\pvettore{y\\x^2} \oplus f(a) = f\left(\pvettore{y\\x^2} \oplus a\right) = f\pvettore{y^1\\x^2}$$

Case 3): General case. Iterating the proof of case 2 we can prove the lemma for $y = y^1 - n(x^1-y^1)$. Then by case 1 we can extend the result to every $y$.
\end{proof}

\begin{proposition}
Let $M$ be a $\tro$-semimodule and let $f:\tro\pro^1 \freccia \pro(M)$ be a tropical projective map. If $f$ is not injective there are two points $x \prec y \in \tro\pro^1$ such that $f(x) = f(y) = p \in \pro(M)$. Then either $\forall z \prec y: f(z) = p$ or $\forall z \succ x: f(z) = p$.
\end{proposition}

\begin{proof}
The map $f$ is associated with a map $\bar{f}: \tro^2 \freccia M$. There exists lifts $\bar{x}, \bar{y} \in \tro^2$ such that $\bar{f}(\bar{x}) = \bar{f}(\bar{y}) = \bar{p}$. Now:

Case 1) If $\bar{x} \preceq \bar{y}$ then one of their coordinates is equal. Else there is a scalar $\lambda < 1_\tro$ such that $x \preceq \lambda \odot y$, and $\bar{p} \preceq \lambda \odot \bar{p}$, a contradiction. Then we can apply the previous lemma, and we have that $\forall z \prec y: f(z) = p$.

Case 2) If $\bar{y} \preceq \bar{x}$ as before we have $\forall z \succ x: f(z) = p$.

Case 3) If they are not comparable, then both are minor than their sum, $\bar{x} \oplus \bar{y}$, and $\bar{f}(\bar{x} \oplus \bar{y}) = p$. Then, by previous cases we have that $\forall z \in \tro\pro^1: f(z) = p$.
\end{proof}

\begin{corollary}    \label{corol:segments}
Let $M$ be a $\tro$-semimodule and let $f:\tro\pro^1 \freccia \pro(M)$ be a tropical projective map. The sets $f^{-1}(f(0_\tro))$ and $f^{-1}(f(\infty_\tro))$ are, respectively, an initial and a final segment for the order of $\tro\pro^1$. If $M$ is separated, then these segments are closed segments. On the complement of their union the map is injective.
\end{corollary}

Suppose that $M$ is a separated $\tro$-semimodule, $\bar{f}:\tro^n \freccia M$ is a linear map and $f:\tro\pro^{n-1} \freccia \pro(M)$ is the induced projective map. As usual we denote by $e_1, \dots, e_n$ the points of the canonical basis of $\tro^n$, and we pose $v_i = \bar{f}(e_i) \in M$. We want to describe the set $V_i = f^{-1}(\pi(v_i))$. It is enough to describe $V_1$. As $\Span_\tro(e_j,e_1)$ is isomorphic to $\tro^2$, we know that $S_j = V_1 \cap \pi(\Span_\tro(e_j,e_1))$ is a closed initial segment of $\pi(\Span_\tro(e_j,e_1))$, with extreme point $\pi(w_j)$. We can suppose that $w_j = a_j e_j + e_1$.

\begin{lemma}   \label{lemma:finitely gen}
The set $V_1$ is 
$$ \pi(\{e_1 \oplus b_2\odot e_2 \oplus \dots \odot b_n\oplus e_n \ |\ b_i \leq a_i \}) $$
Hence there is a point $h_1 = e_1 \oplus a_2\odot e_2 \oplus \dots \odot a_n\oplus e_n$ such that $\pi(h_1)$ is an extremal point of $V_1$.

The restriction of $\bar{f}$ to the submodule $\Span_\tro(h_i,h_j)$ is injective.
\end{lemma}

\subsection{Definition of the metric}

Any convex subset $C$ of a real projective space $\erre\pro^n$ has a well defined metric, the Hilbert metric. This metric is defined by using cross-ratios: if $x,y \in C$, the projective line through $x$ and $y$ intersects $\partial C$ in two points $a,b$. The distance is then defined as $d(x,y) = \frac{1}{2}\log[a,x,y,b]$ (order chosen such that $\overline{ax} \cap \overline{yb} = \emptyset$). If $C,D$ are convex subsets of $\erre\pro^n$ and if $f:C \freccia D$ is the restriction of a projective map, then $d(f(x),f(y)) \leq d(x,y)$. In particular any projective isomorphism $f:C \freccia C$ is an isometry. Moreover this metric has straight lines as geodesics. See \cite{Ki01} for a reference on the Hilbert metric in relation with projective structures.

We can give an analogous definition for separated tropical projective spaces over $\tro_\erre$. In the following we will assume $\Lambda = \erre$ and $\tro = \tro_\erre$. If $M$ is a separated $\tro$-module there is a canonical way for defining a distance $d:\pro(M)\times \pro(M) \freccia \erre \cup \{+\infty\}$. This distance differs from ordinary distances as it can take the value $+\infty$, but has the other properties of a distance (non degeneracy, symmetry, triangular inequality). If $f:\pro(M) \freccia \pro(N)$ is a projective map, then $d(f(x),f(y)) \leq d(x,y)$, and if $S \subset M$ is such that $f_{|S}$ is injective, then $f_{|S}$ is an isometry.

This metric can be defined searching for a tropical analogue of the cross ratio. In $\erre\pro^1$ the cross ratio can be defined by the identity $[0,1,z,\infty] = z$ and the condition of being a projective invariant. Or equivalently if $A$ is the (unique) projective map satisfying $A(0) = a, A(1) = b, A(\infty) = d$, then $[a,b,c,d] = A^{-1}(c)$. In this form the definition can be transposed to the tropical case. 

Let $\tro$ be a tropical semifield and let $a = [a^1:a^2], b = [b^1:b^2], c = [c^1:c^2], d = [d^1:d^2] \in \tro\pro^1 = \pro(\tro^2)$ be points such that $a \prec b \prec c \prec d$. There is a unique tropical projective map $A$ satisfying $A(0_\tro) = a, A(1_\tro) = b, A(\infty_\tro) = d$. This map is described by the matrix 
$$
\begin{gmatrice}
a^2\odot b^1 \odot d^1    &   a^1\odot b^2 \odot d^1\\
a^2\odot b^1 \odot d^2    &   a^2\odot b^2 \odot d^1
\end{gmatrice}
$$ 

Given an $x \in \tro, x \geq 1_\tro$, we have that 
$$A([x:1_\tro]) = \left\{
\begin{array}{ll}
[b^1+x:b^2] &  \mbox{ if } x < (d^1-d^2)-(b^1-b^2)\\
d           &  \mbox{else}
\end{array}
\right.$$

The point $A^{-1}(c)$ is then $[(c^1-c^2)-(b^1-b^2):1]$. Then we can define this point of $\tro\pro^1$ as the cross-ratio of $[a,b,c,d]$. This value depends only on the central points $b,c$, and it is invariant by every tropical projective map $B:\tro\pro^1 \freccia \tro\pro^1$ that is injective on the interval $[b,c]$. 

Consider a tropical projective map $B:\tro\pro^1 \freccia \tro\pro^1$ such that $B(0_\tro) = b$ and $B(\infty_\tro) = c$. This map is described by a matrix of the form:
$$
\begin{gmatrice}
\mu \odot c^1   &   \lambda \odot b^1\\
\mu \odot c^2   &   \lambda \odot b^2
\end{gmatrice}
$$
The inverse images $B^{-1}(b)$ and $B^{-1}(c)$ are, respectively, an initial segment and a final segment of $\tro\pro^1$ with reference to the order $\preceq$ of $\tro\pro^1$. This segments have an extremal point, $b_0$ and $c_0$ respectively. The restriction $B_{|[b_0,c_0]}:[b_0,c_0] \freccia [b,c]$ is a projective isomorphism, hence $(c^1-c^2)-(b^1-b^2) = (c_0^1-c_0^2)-(b_0^1-b_0^2)$.

When we define the Hilbert metric we don't need to take the logarithms, as coordinates in tropical geometry already are in logarithmic scale. Hence the Hilbert metric on $\tro\pro^1$ is simply the Euclidean metric: 
$$ d(x,y) = |(x^1-x^2) - (y^1-y^2)| $$

This definition can be extended to every separated tropical projective space $\pro(M)$. If $a,b \in \pro(M)$, we can choose two lifts $\bar{a}, \bar{b} \in M$. Then there is a unique linear map $\bar{f}:\tro^2 \freccia M$ such that $f(e_1) = \bar{b}, f(e_2) = \bar{a}$. The induced projective map $f:\tro\pro^1 \freccia \pro(M)$ sends $0_\tro$ in $a$ and $\infty_\tro$ in $b$. By corollary \ref{corol:segments} the sets $f^{-1}(a)$ and $f^{-1}(b)$ are closed segments, with extremal points $a_0, b_0$. We can define the distance as $d(a,b) = d(a_0, b_0)$.
It is easy to verify that this definition does not depend on the choice of the lifts $\bar{a}, \bar{b}$. Now we have to verify the triangular inequality, but it is more comfortable to give an example first.

For the projective spaces associated with the free modules we can calculate explicitly this distance. It is a well known distance, the Hilbert metric on the simplex in logarithmic coordinates.

\begin{proposition}
Let $x, y \in \tro\pro^{n-1}$. Then, for all lifts $\bar{x}, \bar{y} \in \tro^n$: 
$$ d(x,y) = \left( \bigoplus_{i=1}^n \bar{x}^i \oslash \bar{y}^i \right) \odot \left( \bigoplus_{i=1}^n \bar{y}^i \oslash \bar{x}^i \right) = \max_{i=1}^n (\bar{x}^i-\bar{y}^i) + \max_{i=1}^n(\bar{y}^i - \bar{x}^i)$$
\end{proposition}

\begin{proof}
The map $\bar{f}$ as above is defined in this case by the $2 \times n$ matrix:
$$\begin{gmatrice}
y^1   &   x^1\\
\vdots &  \vdots \\
y^n  &  x^n  
\end{gmatrice}$$
Then for all $h \in \tro$, 
$${\left(\bar{f}\pvettore{h\\1_\tro}\right)}^i = \max(y_i\odot h, x^i)$$
This is equal to $x$ if $\forall i: h \leq x^i - y^i$, i.e. if $h \leq \min_i (x^i - y^i)$.
As before, for all $k \in \tro$,
$${\left(\bar{f}\pvettore{1_\tro\\k}\right)}^i = \max(y_i, x^i\odot k)$$
This is equal to $y$ if $\forall i: k \leq y^i - x^i$, i.e. if $k \leq \min_i (y^i - x^i)$
Then, by definition
$$d(x,y) = |\min_i (x^i - y^i) + \min_i (y^i - x^i)| $$
By changing signs inside the absolute value, we have the thesis.
\end{proof}

From this explicit computation we can deduce easily that the triangular inequality holds for the distance we have defined in $\tro\pro^{n-1}$, and that the topology induced by this distance on $\tro\pro^{n-1}$ is the quotient of the product topology on $\tro^n$.

Once we know that the triangular inequality holds for $\tro\pro^{n-1}$, we can use this fact to prove it for all separated tropical projective spaces.

\begin{proposition}
Let $M$ be a separated $\tro$-semimodule. Then the function $d:\pro(M) \times \pro(M) \freccia \erre_{\geq 0} \cup \{\infty\}$ satisfy:
$$\forall x,y,z \in \tro\pro^{n-1}: d(x,y) \leq d(x,z)+ d(z,y) $$
\end{proposition}

\begin{proof}
Fix lifts $\bar{x}, \bar{y}, \bar{z} \in M$. We can construct a map $f:\tro^3 \freccia M$ such that $f(e_1) = x, f(e_2) = y, f(e_3) = z$. By lemma \ref{lemma:finitely gen} there exist points $h_1,h_2,h_3\in \tro^3$ such that $f$ is injective over $\Span_\tro(h_i,h_j)$. Then $d(\pi(h_i),\pi(h_j)) = d(\pi(f_i),\pi(f_j))$. As the triangular inequality holds in $\tro\pro^2$, then it holds for $x,y,z$.
\end{proof}

The metric we have defined for separated tropical projective spaces can achieve the value $+\infty$. Given a $\tro$-semimodule $M$ we can define the following equivalence relation on $M \setminus \{0\}$:
$$x \sim y \Leftrightarrow d(\pi(x),\pi(y)) < +\infty$$
The union of $\{0\}$ with one of these equivalence classes is again a $\tro$-semimodule, and their projective quotients are tropical projective spaces with an ordinary (i.e. finite) metric.

For example, if $M = \tro^n$, the equivalence class of the point $(1_\tro, \dots, 1_\tro)$ is the $\tro$-semimodule $F\tro^n$, and its associated projective space is $F\tro\pro^{n-1}$, a tropical projective space in which the metric is finite.

For the $\tro$-semimodule $U^n$ an equivalence class is $FU^n$, and its associated projective space is  $FP^{n-1}$, a tropical projective space in which the metric is finite.

We can calculate more explicitly the metric for $FP^{n-1}$. Let $x,y\in FP^{n-1}$ and let $\bar{x}, \bar{y}\in U^n$ be their lifts. By lemma \ref{lemma:common basis} there exists a basis $\ecors = (e_1, \dots, e_n)$ of $\bar{x}$ such that $a_1 e_1, \dots, a_n e_n$ is a basis of $\bar{y}$. In the tropical chart $\varphi_\ecors$, the point $\bar{x}$ has coordinates $(1_\tro,\dots,1_\tro)$, while the point $\bar{y}$ has coordinates $(\tau(a_1),\dots,\tau(a_n))$. Hence 
$$d(x,y) = \max_i(\tau(a_i))-\min_i(\tau(a_i)) $$

\subsection{Homotopy properties}

In this section we will show that every separated tropical projective space with a finite metric is contractible.

If $(X,d)$ is a metric space, we denote by $C^0([0,1],X)$ the space of continuous curves in $X$, with the metric defined by 
$$d(\gamma,\gamma') = \max_{t \in [0,1]} d(\gamma(t),\gamma'(t))$$ 
Note that the following pairing is continuous
$$C^0([0,1],X) \times [0,1] \ni (\gamma,t) \freccia \gamma(t) \in X$$

\begin{lemma}   \label{lemma:contractible}
Let $(X,d)$ be a metric space and suppose we can construct a continuous map:
$$C:X\times X \ni (x,y) \freccia C_{x,y} \in C^0([0,1],X) $$
such that
\begin{enumerate}
\item $C_{x,y}(0) = x$ and $C_{x,y}(1) = y$ 
\item $C_{x,x}$ is a constant curve.
\end{enumerate}
Then $X$ is contractible.
\end{lemma}

\begin{proof}
We can construct a retraction $ H : X \times [0,1] \freccia X $ retracting $X$ on one of its points $\{\bar{x}\}$ as
$$H(y,t) = C_{y,\bar{x}}(t)$$
By definition of $C$ we have that $H(y,0) = y$ and $H(y,1) = \bar{x}$, and $H$ is continuous as it is a composition of continuous functions.
\end{proof}

\begin{lemma}   \label{lemma:continuity}
Let $x,y,a,b \in \tro^n$ and let $\phi_{x,a}$ and $\phi_{y,b}$ be, respectively, the linear maps $\tro^2 \freccia \tro^n$ defined by the matrices:
$$ \phi_{x,a} = 
\begin{gmatrice}
x^1  &  a^1\\
\vdots & \vdots\\
x^n  &  a^n
\end{gmatrice}, \hspace{2cm}
\phi_{y,b} = 
\begin{gmatrice}
y^1  &  b^1\\
\vdots & \vdots\\
y^n  &  b^n
\end{gmatrice}$$
Then $$\forall v \in \tro^2: d(\pi(\phi_{x,a}(v)), \pi(\phi_{y,b}(v))) \leq \max( d(\pi(x),\pi(y)), d(\pi(a),\pi(b)))$$
\end{lemma}

\begin{proof}
Without loss of generality we can suppose that $v = (t,1_\tro)$, so that $(\phi_{x,a}(v))^i = \max(x^i+t,a^i)$. Then
$$d(\pi(\phi_{x,a}(v)), \pi(\phi_{y,b}(v))) = $$ $$=\max_i(\max(x^i+t,a^i) - \max(y^i+t,b^i)) + \max_i(\max(y^i+t,b^i)) - \max(x^i+t,a^i)$$
It is easy to check that $\max(x^i+t,a^i) - \max(y^i+t,b^i)) \leq \max(x^i-y^i,a^i-b^i)$ by analyzing the four cases.
\end{proof}

\begin{proposition}
For every separated $\tro$-module $M$, its associated projective space $\pro(M)$ is contractible with reference to the topology induced by the canonical metric.
\end{proposition}

\begin{proof}
We have to construct a map $C$ as in lemma \ref{lemma:contractible}. We will use tropical segments, rescaling their parametrization to the interval $[0,1]$. If $x,y \in \pro(M)$, we take lifts $\bar{x},\bar{y} \in M$ and the map $\bar{f}:\tro^2 \freccia M$ s.t. $\bar{f}(e_1) = x, \bar{f}(e_2) = y$. As usual $f: \tro\pro^1 \freccia \pro(M)$ is the induced map. By corollary \ref{corol:segments} the sets $f^{-1}(x)$ and $f^{-1}(y)$ are closed segments, with extremal points $x_0, y_0$, hence $f$ restricted to the interval $[x_0,y_0]$ is a curve joining $x$ and $y$. Let $\phi$ be the affine map from the interval $[x_0,y_0]$ to the interval $[0,1]$. We define $C_{x,y}$ as the reparametrization of $f$ by $\phi$. Properties 1 and 2 of the lemma \ref{lemma:contractible} holds for $C$. To prove 3 we can show that:
$$\forall x,y,z,w \in \pro(M): \forall t \in [0,1]: d(C_{x,y}(t),C_{z,w}(t)) \leq 3 \max(d(x,z),d(y,w))$$
To do this we take lifts $\bar{x},\bar{y},\bar{z},\bar{w} \in M$, and a map $\bar{f}: \tro^4 \freccia M$ s.t. $f(e_1) = x, f(e_2) = y, f(e_3)=z, f(e_4)=w$. By lemma \ref{lemma:finitely gen} there exist points $h_1,h_2,h_3,h_4\in \tro^4$ such that $f$ is injective over $\Span_\tro(h_i,h_j)$. Then $d(\pi(h_i),\pi(h_j)) = d(\pi(f_i),\pi(f_j))$. Moreover $f$ is $1$-Lipschitz on $\pi(\Span_\tro(h_1,\dots,h_4))$, hence our property on $M$ follows from the same property on $\tro^4$, and this follows from lemma \ref{lemma:continuity}.
\end{proof}

\section{Tropicalization of group representations}  \label{sec:Tropicalization of group actions}

\subsection{Length spectra}

Let $\Gamma$ be a group and $\rho:\Gamma \freccia GL_{n+1}(\effe)$ be a representation of $\Gamma$ in the general linear group of a non-archimedean field with surjective real valuation.

Let $\effe$ be a non-archimedean field with surjective real valuation. The group $GL_{n+1}(\effe)$ acts by linear maps on the tropical modules $U^{n+1}(\effe)$ and $FU^{n+1}(\effe)$, and by tropical projective maps on the tropical projective spaces $P^n(\effe)$ and $FP^n(\effe)$. The representation $\rho$ defines an action of $\Gamma$ on $FP^n(\effe)$.

For every matrix $A \in GL_{n+1}(\effe)$, we can define the \nuovo{translation length} of $A$ as:
$$l(A) = \inf_{x \in FP^n(\effe)} d(x,Ax)$$

\begin{proposition}   \label{prop:dist}
Let $x \in FP^n$, and $L \subset V$ be a lift of $x$ in $FU^{n+1}$. We denote by $e_1, \dots, e_{n+1}$ a basis of $L$, and by $\tilde{A}$ the matrix corresponding to $A$ in this basis. Then
$$d(x,A(x)) =  \max_{i,j} \tau({(\tilde{A})}^i_j) + \max_{i,j} \tau({(\tilde{A}^{-1})}^i_j)$$
\end{proposition}

\begin{proof}
By lemma \ref{lemma:common basis} applied to the $\ocors$-modules $L$ and $A(L)$, there exist a basis $v_1, \dots, v_n$ of $L$ and scalars $\lambda_1, \dots, \lambda_n \in \effe$ such that $\lambda_1 v_1, \dots, \lambda_n v_n$ is a basis of $A(L)$. Then $d(x,Ax) = \max_i(\tau(\lambda_i))-\min_i(\tau(\lambda_i))$. We will denote by $M_1$ the transition matrix from $e_1,\dots,e_n$ to $v_1,\dots,v_n$. As they are bases of the same $\ocors$-module, $M_1$ is in $GL_n(\ocors)$. We will denote by $M_2$ the transition matrix from $\lambda_1 v_1, \dots, \lambda_n v_n$ to $A(e_1), \dots, A(e_n)$, and it is again in $GL_n(\ocors)$. Let $\Delta$ be the diagonal matrix:
$$\Delta = \begin{gmatrice}
\lambda_1 & & \\
  & \ddots & \\
  &    &  \lambda_n
\end{gmatrice}$$
Then the following relations hold:
$$\tilde{A}= M_2 \Delta M_1 \ \ \ \ \Delta = M_2^{-1} \tilde{A} M_1^{-1}$$ 
$$\tilde{A}^{-1}= M_1^{-1} \Delta^{-1} M_2^{-1} \ \ \ \ \Delta^{-1} = M_1 \tilde{A}^{-1} M_2$$
Hence:
$$\max_i(\tau(\lambda_i)) = \max_{i,j} \tau({(\tilde{A})}^i_j)$$
In the same way we have:
$$-\min_i(\tau(\lambda_i)) = \max_{i,j} \tau({(\tilde{A}^{-1})}^i_j)$$
\end{proof}

The case $n=1$ has been studied in \cite{MS84}. If $A \in SL_2(\effe)$, we have $l(A) = 2\max(0,\tau(\trace(A)))$ (see \cite[prop. II.3.15]{MS84}). In the following we give an extension of this result for generic $n$.

Let $\effe$ be a non-archimedean real closed field of finite rank extending $\erre$, with a surjective real valuation $\overline{v}: \effe^* \freccia \erre$ such that the valuation ring is convex. The field $\cappa = \effe[i]$ is an algebraically closed field extending $\ci$, with an extended valuation $\overline{v}:\cappa^* \freccia \erre$. We will use the notation $\tau = -\overline{v}$. We will also use the complex norm $|\cdot|:\cappa \freccia \effe_{\geq 0}$ defined by $|a+bi| = \sqrt{a^2 + b^2}$ and the conjugation $\overline{a+bi} = a-bi$.

If $A \in GL_{n+1}(\cappa)$, we denote by $\lambda_1, \dots, \lambda_{n+1}$ its eigenvalues, ordered such that $|\lambda_i| \geq |\lambda_{i+1}|$. We will denote $r(A) = |\lambda_1|$, the \nuovo{spectral radius} of $A$.

Note that the function 
$$|\cdot|: M_n(\cappa) \ni A \freccia \in \max_{i,j}|A^i_j| \in \effe_{\geq 0}$$
is a consistent norm on $M_n(\cappa)$, hence, by the spectral radius theorem, we have $r(A) \leq |A|$.

\begin{proposition}   \label{prop:spectrum_complex}
Suppose the field $\cappa$ is as above. Then a matrix $A \in GL_{n+1}(\cappa)$ acts on $FP^n(\cappa)$. Then the $\inf$ in the definition of $l(A)$ is a minimum, and it is equal to
$$l(A) = \tau\left(\left|\frac{\lambda_1}{\lambda_{n+1}}\right|\right)$$
\end{proposition}

\begin{proof}
By proposition \ref{prop:dist} we have that for every $x \in FP^n(\cappa)$
$$d(x,A(x)) \geq \tau(r(A)) + \tau(r(A^{-1}))$$
Hence
$$l(A) \geq  \tau(r(A)) + \tau(r(A^{-1}))$$
or, in other words,
$$l(A) \geq \tau\left(\left|\frac{\lambda_1}{\lambda_{n+1}}\right|\right)$$
We only need to show that the lower bound of previous corollary is actually achieved. The Jordan form of $A$ is
$$
\begin{gmatrice}
\lambda_1 & *   \\
          & \lambda_2  & \ddots  \\
          &            & \ddots & * \\
          &            &        & \lambda_{n+1} \\
\end{gmatrice}
$$
where the entries marked by $*$ are $0$ or $1$. Let $v_1,\dots,v_{n+1}$ be a Jordan basis, and let $L=\Span_\ocors(v_1,\dots,v_n) \in U^{n+1}$. By proposition \ref{prop:dist}
$$d(\pi(L),A\pi(L)) = \tau\left(\left|\frac{\lambda_1}{\lambda_{n+1}}\right|\right)$$
\end{proof}

Now suppose that $A \in GL_n(\effe)$, with $\effe$ a non-archimedean real closed field as above. Hence $A$ acts on $FP^n(\effe)$, and now we want to study the translation length of $A$ over $FP^n(\effe)$. As before, we denote by $\lambda_1, \dots, \lambda_{n+1} \in \cappa$ its eigenvalues, ordered such that $|\lambda_i| \geq |\lambda_{i+1}|$.

\begin{proposition}  \label{prop:spectrum_real}
Suppose that $\effe$ is as above, and that $A \in GL_{n+1}(\effe)$. We consider the translation length $l(A)$ with respect to the action of $A$ on $FP^n(\effe)$. Then the $\inf$ in the definition of $l(A)$ is a minimum, and it is equal to
$$l(A) = \tau\left(\left|\frac{\lambda_1}{\lambda_{n+1}}\right|\right)$$
\end{proposition}

\begin{proof}
As $FP^n(\effe) \subset FP^n(\cappa)$, by proposition \ref{prop:spectrum_complex} we have the inequality
$$l(A) \geq \tau\left(\left|\frac{\lambda_1}{\lambda_{n+1}}\right|\right)$$
To prove that this lower bound is achieved, we will choose a suitable basis, as above. Consider the decomposition into sum of generalized eigenspaces
$$\cappa^{n+1} = \sum_{i = 1}^n \ker({(A-\lambda_i \ident)}^{n+1}) $$
For every $\lambda_i \in \effe$, the generalized eigenspace $\ker({(A-\lambda_i \ident)}^{n+1})$ has a basis of generalized eigenvectors in $\effe^{n+1}$. If $\lambda_i \in \cappa \setminus \effe$, then $\overline{\lambda_i}$ is an eigenvalue, and if $v_1, \dots, v_s$ is a basis of generalized eigenvectors of $\ker({(A-\lambda_i \ident)}^{n+1})$, then $\overline{v_1}, \dots, \overline{v_s}$ is a basis of generalized eigenvectors of $\ker({(A-\overline{\lambda_i} \ident)}^{n+1})$. The vectors $v_i + \overline{v_i}$ and $v_i - \overline{v_i}$ are in $\effe^{n+1}$, and they form a basis of $\ker({(A-\lambda_i \ident)}^{n+1}) + \ker({(A-\overline{\lambda_i} \ident)}^{n+1})$. In this way we have constructed a basis $v_1,\dots,v_{n+1}$  of $\effe^{n+1}$ such that $|A| = |\lambda_1|$ and $|A^{-1}| = |\lambda_{n+1}|$. Let $L=\Span_\ocors(v_1,\dots,v_n) \in U^{n+1}$, then, by proposition \ref{prop:dist}
$$d(\pi(L),A\pi(L)) = \tau\left(\left|\frac{\lambda_1}{\lambda_{n+1}}\right|\right)$$
\end{proof}

\subsection{Boundary points}

Here we give a geometric interpretation to the points of the boundaries of the spaces of convex projective structures. Let $M$ be a closed $n$-manifold such that the fundamental group $\pi_1(M)$ has trivial virtual center, it is Gromov hyperbolic, and it is torsion free (note that every closed hyperbolic $n$-manifold whose fundamental group is torsion-free satisfies the hypotheses). In \cite[subsec. 6.4]{A2}, we considered the family $\gfamil = {\{e_\gamma\}}_{\gamma \in \pi_1(M)}$, and we constructed a compactification of $\mathcal{T}^{c}_{\erre\pro^n}(M)$:
$$\overline{\mathcal{T}^{c}_{\erre\pro^n}(M)}_{\gfamil} = \mathcal{T}^{c}_{\erre\pro^n}(M) \cup \partial_\gfamil \mathcal{T}^{c}_{\erre\pro^n}(M)$$
The cone over the boundary $C(\partial_\gfamil \mathcal{T}^{c}_{\erre\pro^n}(M))$ can be identified with a subset of $\erre^{\gfamil} = \erre^{\pi_1(M)}$.

Every action of $\pi_1(M)$ on a tropical projective space $FP^n(\effe)$ has a well defined length spectrum $(l(\gamma))_{\gamma \in \pi_1(M)} \in \erre^{\pi_1(M)}$.

\begin{theorem}
Let $\effe = \erre((t^{\erre^r}))$, where $r$ is the dimension of $\mathcal{T}^{c}_{\erre\pro^n}(M)$ (see the definition in \cite[subsec. 3.3]{A2}). The points of $C(\partial_\gfamil \mathcal{T}^{c}_{\erre\pro^n}(M))$ are length spectra of actions of the fundamental group $\pi_1(M)$ on the tropical projective space $FP^n(\effe)$.
\end{theorem}

\begin{proof}
The semi-algebraic set $\mathcal{T}^{c}_{\erre\pro^n}(M)$ has an extension to the field $\effe$, that we will denote by $\overline{\mathcal{T}^{c}_{\erre\pro^n}(M)} \subset \overline{\charat}(\pi_1(M), SL_{n+1}(\effe))$. Every element of $\overline{\mathcal{T}^{c}_{\erre\pro^n}(M)}$ is a conjugacy class of a representation $\rho: \pi_1(M) \freccia SL_{n+1}(\effe)$.

Let $x \in C(\partial_\gfamil \mathcal{T}^{c}_{\erre\pro^n}(M)) \subset \erre^{\gfamil}$. As we said in \cite[subsec. 3.3]{A2}), there exists a representation $\rho \in \overline{\mathcal{T}^{c}_{\erre\pro^n}(M)}$ such that for every $\gamma \in \pi_1(M)$, the matrix $\rho(\gamma)$ satisfies $x_{e_\gamma} = \tau\left(\left|\frac{\lambda_1}{\lambda_{n+1}}\right|\right)$.

Consider the action of $\pi_1(M)$ on $FP^{n}(\effe)$ induced by the representation $\rho$. By proposition \ref{prop:spectrum_real}, the translation length of an element $\gamma$ is
$$l(\rho(\gamma)) = \tau\left(\left|\frac{\lambda_1}{\lambda_{n+1}}\right|\right)$$
\end{proof}

This result is an extension of the interpretation of the boundary points of the Teichm\"uller spaces given by Morgan and Shalen in \cite{MS84}. Here we review their result in our language. Let $\overline{S} = \Sigma_g^k$, a surface of genus $g$ with $k \geq 0$ boundary components and such that $\chi(\overline{S}) < 0$, and let $S$ be the interior part of $\overline{S}$. In \cite[subsec. 5.2]{A2} we considered the family $\gfamil = {\{J_\gamma\}}_{\gamma \in \pi_1(S)}$, and we constructed a compactification of $\mathcal{T}^{cf}_{\iper^2}(S)$:
$$\overline{\mathcal{T}^{cf}_{\iper^2}(S)}_{\gfamil} = \mathcal{T}^{cf}_{\iper^2}(S) \cup \partial_\gfamil \mathcal{T}^{cf}_{\iper^2}(S)$$
The cone over the boundary $C(\partial_\gfamil \mathcal{T}^{cf}_{\iper^2}(S))$ can be identified with a subset of $\erre^{\gfamil} = \erre^{\pi_1(S)}$.

Every action of $\pi_1(S)$ on a tropical projective space $FP^1(\effe)$ has a well defined length spectrum $(l(\gamma))_{\gamma \in \pi_1(S)} \in \erre^{\pi_1(S)}$.

\begin{theorem}
Let $\effe = H(\overline{\erre}^\erre)$, the Hardy field as in \cite[subsec. 4.1]{A1}. The points of $C(\partial_\gfamil \mathcal{T}^{cf}_{\iper^2}(S))$ are length spectra of actions of the fundamental group $\pi_1(S)$ on the tropical projective space $FP^1(\effe)$.
\end{theorem}

\begin{proof}
The semi-algebraic set $\mathcal{T}^{cf}_{\iper^2}(M)$ has an extension to the field $\effe$, denoted by $\overline{\mathcal{T}^{cf}_{\iper^2}(S)} \subset \overline{\charat}(\pi_1(S), SL_2(\effe))$. Given a representation $\rho \in \overline{\mathcal{T}^{cf}_{\iper^2}(S)}$, if $\gamma \in \pi_1(S)$, the matrix $\rho(\gamma)$ has $|\trace(\rho(\gamma))| \geq 2$.

Let $x \in C(\partial_\gfamil \mathcal{T}^{cf}_{\iper^2}(S))$. By \cite[prop. 3.6]{A2}, there exists a representation $\rho \in \overline{\mathcal{T}^{cf}_{\iper^2}(S)}$ such that for every $\gamma \in \pi_1(S)$, the matrix $\rho(\gamma) \in SL_2(\effe)$ satisfies $x_{J_\gamma} = \tau\left(\trace(\rho(\gamma))\right)$.

Consider the action of $\pi_1(S)$ on $FP^{n}(\effe)$ induced by the representation $\rho$. By \cite[prop. II.3.15]{MS84}, the translation length of every $\gamma \in \pi_1(M)$ is $l(\rho(\gamma)) = 2\max(0,\tau(\trace(A)))$. As $|\trace(\rho(\gamma))| \geq 2$, we have
$l(\rho(\gamma)) = 2 \tau(\trace(A))$
\end{proof}

\subsection{The equivariant map}

These actions of $\pi_1(M)$ on the tropical projective spaces $FP^{n}$ should correspond to some kind of dual structure on $M$.

Suppose that $M$ is an $n$-manifold such that, if $n > 2$, $\pi_2(M) = \dots = \pi_{n-1}(M) = 0$. For example every manifold whose universal covering is $\erre^n$ satisfy this hypothesis, in particular every manifold admitting a convex projective structure. We will denote by $p:\tilde{M} \freccia M$ the universal covering of $M$. Then suppose that $Z$ is a simply connected topological space with an action of $\pi_1(M)$. It is always possible to construct an equivariant map:

\begin{theorem}           \label{teo:equivariant map}
There exists a map $f:\tilde{M} \freccia Z$ that is equivariant for the action of $\pi_1(M)$, i.e.
$$\forall x \in \tilde{M}: \forall \gamma \in \pi_1(M): \gamma(f(x))=f(\gamma(x)) $$
\end{theorem}

\begin{proof}
The group $\pi_1(M)$ acts diagonally on the space $\tilde{M} \times Z$: 
$$\gamma(x,z) = (\gamma(x),\gamma(z))$$ 
This action is free and proper, $\tilde{M} \times Z$ is simply connected, hence 
$$P:\tilde{M} \times Z \freccia \mathcal{K} = (\tilde{M} \times Z) / \pi_1(M)$$ 
is a universal cover, and $\pi_1(\mathcal{K}) = \pi_1(M)$.

As $M$ is a manifold it is homeomorphic to a CW-complex of dimension $n$ with only one $0$-cell. Hence the hypothesis that $\pi_2(M) = \dots = \pi_{n-1}(M) = 0$ implies that the isomorphism $\pi_1(M) \freccia \pi_1(\mathcal{K})$ is induced by a map $\psi:M \freccia \mathcal{K}$, well defined up to homotopy. 

As $\tilde{M}$ is simply connected, we can lift the map $\phi = \psi \circ p: \tilde{M} \freccia \mathcal{K}$ to a map $\tilde{\phi}: \tilde{M}\freccia \tilde{M} \times Z$ such that $P \circ \tilde{\phi} = \phi$. The equivariant map $f$ we are searching for is the composition of $\tilde{\phi}$ with the projection on $Z$. We have to check that it is equivariant, and to show this we will prove that $\tilde{\phi}$ is equivariant. We need to prove that:
$$\forall y \in \tilde{M}: \forall \gamma \in \pi_1(M): \gamma(\tilde{\phi}(y)) = \tilde{\phi}(\gamma(y)) $$
Fix an $y \in \tilde{M}$ and a $\gamma \in \pi_1(M)$. Let $x_0 = p(y) = p(\gamma(y)) \in M$ and let $x_1 = \psi(x_0) = P(\tilde{\phi}(y)) = P(\tilde{\phi}(\gamma(y)))$ (as $\tilde{\phi}$ is a lift of $\psi:M \freccia \mathcal{K}$).

Now we identify $\pi_1(M)$ with the based fundamental groups $\pi_1(M,x_0)$ and $\pi_1(\mathcal{K},x_1)$. By the definition of $\psi$, with this identification, the isomorphism $\psi_*:\pi_1(M,x_0) \freccia \pi_1(\mathcal{K},x_1)$ is the identity, hence $\psi_*(\gamma)=\gamma$.

Consider the lift $\tilde{\gamma}$ of the path $\gamma$ in $\tilde{M}$ starting from the point $y$. The other extreme of $\tilde{\gamma}$ is the point $\gamma(y)$. The same way the lift $\widetilde{\psi_*(\gamma)}$ of the path $\psi_*(\gamma)$ in $\tilde{M} \times \mathcal{K}$ starting from the point $\tilde{\phi}(y)$ is the image $\tilde{\phi}(\tilde{\gamma})$, hence the other extreme of this path is the point $\tilde{\phi}(\gamma(y))$. This is precisely the definition of $\gamma(\tilde{\phi}(y))$.
\end{proof}

Suppose that $M$ is as above, and that we have an action of $\pi_1(M)$ on the tropical projective space $P^m$. As $P^m$ is a contractible space there is a $\pi_1(M)$-equivariant map
$$f:\widetilde{M} \freccia P^n$$
An interesting open problem is to understand the dual structure this equivariant map induces on $M$.

The case where $M$ is an hyperbolic surface and $m=1$ has been studied by Morgan and Shalen in \cite{MS88} and it is well understood: $P^1$ is a real tree and the equivariant map induces a measured lamination on $M$, that is dual to the action.

This work can possibly lead to the discovery of analogous structures for the general case. For example an action of $\pi_1(M)$ on $P^m$ induces a degenerate metric on $M$, and this metric can be used to associate a length with each curve. Anyway it is not clear up to now how to classify these induced structures. This is closely related to a problem raised by J. Roberts (see \cite[problem 12.19]{Oh01}): how to extend the theory of measured laminations to higher rank groups, such as, for example, $SL_n(\erre)$.

%
%


\begin{thebibliography}{CGQ04}
	\bibitem[A1]{A1} D. Alessandrini, \emph{Logarithmic limit sets of real semi-algebraic sets}, preprint on arXiv:0707.0845 v2, submitted.
	\bibitem[A2]{A2} D. Alessandrini, \emph{A compactification for the spaces of convex projective structures on manifolds}, preprint on arXiv:0801.0165 v1.
	\bibitem[CGQ04]{CGQ04} G. Cohen, S. Gaubert, J.-P. Quadrat, \emph{Duality and separation theorems in idempotent semimodules}, Linear Algebra and its Applications, {\bfseries 379} (2004), 395--422.
	\bibitem[JSY]{JSY} M. Joswig, B. Sturmfels, J. Yu, \emph{Affine Buildings and Tropical Convexity}, preprint on arXiv:0706.1918.
	\bibitem[Ki01]{Ki01} I. Kim, \emph{Rigidity and deformation spaces of strictly convex real projective structures on compact manifolds}, J. Diff. Geometry {\bfseries 58} (2001) 189--218.
	\bibitem[MS84]{MS84} J. W. Morgan, P. B. Shalen, \emph{Valuations, trees, and degeneration of hyperbolic structures, I}, Annals of Math. {\bfseries 120} (1984), 401--476.
	\bibitem[MS88]{MS88} J. W. Morgan, P. B. Shalen, \emph{Degeneration of hyperbolic structures, II: Measured laminations in $3$-manifolds}, Annals of Math. {\bfseries 127} (1988), 403--456.
	\bibitem[MS88']{MS88'} J. W. Morgan, P. B. Shalen, \emph{Degeneration of hyperbolic structures, III: Actions of $3$-manifold groups on trees and Thurston's compactness theorem}, Annals of Math. {\bfseries 127} (1988), 457--519.
	\bibitem[Oh01]{Oh01} T. Ohtsuki, \emph{Problems on invariants of knots and 3-manifolds}, Geometry and Topology Monographs Volume 4: Invariants of knots and 3-manifolds (Kyoto 2001), 377--572.
\end{thebibliography}
\end{document}